\documentclass[aps,prl,reprint]{revtex4-1}

\usepackage{graphicx}
\usepackage{epsfig}
\usepackage{amsmath,amsfonts,amssymb,amsthm,bbm,color}

\theoremstyle{plain}
\newtheorem{thm}{Theorem}
\newtheorem*{thm*}{Theorem}
\newtheorem{lem}{Lemma}
\newtheorem*{lem*}{Lemma}
\newtheorem{cor}{Corollary}
\newtheorem*{cor*}{Corollary}

\newtheorem*{prop*}{Proposition}

\theoremstyle{definition}

\newtheorem{Def}{Definition}

\newcommand{\reftit}{\textit}    
\def\cR{{\mathcal{R}}}

\def\cT{{\mathcal{T}}}

\def\cG{{\mathcal{G}}}

\def\cS{{\mathcal{S}}}

\newcommand\norm[1]{\left\lVert#1\right\rVert}

\def\be{\begin{equation}}
\def\ee{\end{equation}}
\def\ben{\begin{equation*}}
\def\een{\end{equation*}}

\begin{document}

\title{Tokunaga self-similarity arises naturally from time invariance}

\author{Yevgeniy Kovchegov}
\email[]{kovchegy@math.oregonstate.edu}
\affiliation{Department of Mathematics, Oregon State University, Corvallis, OR, 97331-4605, USA}

\author{Ilya Zaliapin}
\email[]{zal@unr.edu}
\affiliation{Department of Mathematics and Statistics, University of Nevada, Reno, NV, 89557-0084,
USA}

\date{\today}

\begin{abstract}
The Tokunaga condition is an algebraic rule that provides a detailed description 
of the branching structure in a self-similar tree.
Despite a solid empirical validation and practical convenience,
the Tokunaga condition lacks a theoretical justification. 
Such a justification is suggested in this work. 
We define a geometric branching processes $\cG(s)$ 
that generates self-similar rooted trees.
The main result establishes the equivalence between the invariance of $\cG(s)$
with respect to a time shift and a one-parametric version of the Tokunaga condition.
In the parameter region where the process satisfies the Tokunaga condition  
(and hence is time invariant), $\cG(s)$ enjoys many of the symmetries observed in a critical binary Galton-Watson 
branching process and reproduce the latter for a particular
parameter value.
\end{abstract}


\maketitle

{\bf
The Tokunaga condition, also called Tokunaga self-similarity, is a particularly
simple parameterization of branching that gives a constructive 
description of a self-similar tree. 
It has originated in hydrology and has been empirically 
rediscovered in a variety of other fields. 
The lack of a theoretical justification behind the Tokunaga condition
raises the question of whether this is an actual physical constraint
or merely a convenient ad hoc approximation that cannot be rejected with available data.
The answer to this question affects theoretical and modeling efforts across
disciplines.
We argue that the Tokunaga condition is an actual physical constraint.
Specifically, we establish a measure-theoretic equivalence between a one-parametric
version of the Tokunaga condition and time invariance of a 
branching process that generates random trees.
The methodology and results of this work are readily applicable to a multitude of questions 
related to scaling laws in trees (Horton law, fractal dimension, etc.). 
We also solve an abstract non-linear problem that can be of general 
interest to the non-linear dynamics community.
}

\section{\label{sec_intro} Introduction}
Tree-shaped fractal formations, from the namesake botanical trees to 
river tributary networks to the systems of canyons 
and mountain crests that define the Earth topography, 
have always fascinated the great minds (e.g., \cite{Richter70}) 
providing inspiration in science, art, and architecture \cite{RS14,Viennot90}.
A quantitative understanding of the branching 
patterns is instrumental in hydrology \cite{RIR01,Tar88,Mar96,BMR99,BM10,Kir00},
geomorphology \cite{DR00,HSP15},
statistical seismology \cite{NTG95,BP04,THR07,HTR08,ZGKW08,G13,ZBZ13,Y13},
statistical physics of fracture \cite{GZNK00,ZKG03a,ZKG03b},  
vascular analysis \cite{Kassab00},
brain studies \cite{Cetal06},
ecology \cite{CLF07}, biology \cite{TPN98}, 
and beyond, encouraging further rigorous treatment.
This study establishes equivalence between the Tokunaga 
self-similarity -- a widely recognized algebraic parameterization
of a self-similar tree -- and a measure-theoretic invariance of a tree 
distribution with respect to a depth shift.

Much of the self-similar tree analysis has originated in the studies of 
river networks; it refers to the rooted trees with no edge lengths.
The celebrated Horton-Strahler ordering scheme \cite{Horton45,Strahler57} 
assigns integer orders to tree vertices and edges, beginning with 
${\sf ord}=1$ at the leaves and incrementing by unity every time 
the two edges of the same order meet (see Def.~\ref{def:HS}, Fig.~\ref{fig:Horton}). 
A sequence of adjacent vertices/edges with the same order is called a 
{\it branch} (Fig.~\ref{fig:Horton}c). 
Let $N_j$ be the number of branches of order $j$ 
and $N_{ij}$, $i<j$, be the number of {\it side branches} of order $\{ij\}$ --
instances when an order-$i$ branch merges 
with an order-$j$ branch in a finite tree $T$. 
The average number of branches of order $i$ in a single branch of order $j$ can 
be traced with the following ratio \cite{KZ16a}
\be
T_{ij}=\frac{{\sf E}[N_{ij}]}{{\sf E}[N_j]}.
\ee
Here ${\sf E}[x]$ is the expected value with respect to a suitable 
probability measure on the examined collection of trees.

Ronald L. Shreve \cite{Shreve69} pioneered the study of side branching in
a {\it topologically random channel network}, which is a uniform distribution
of rooted binary trees with a given number of leaves. 
This model is equivalent to a critical binary Galton-Watson tree, conditioned
on the number of leaves \cite{BWW00,Pitman}.
Recall that a Galton-Watson process describes growth of a population that 
begins with a single progenitor. 
At each discrete time step, each population member disappears
giving birth to a random number $k\ge 0$ of offspring 
according to a distribution $\{p_k\}$. 
A Galton-Watson tree describes a trajectory of this process.
The members are represented by the tree vertices;
the progenitor corresponds to the tree root; the edges connect
parents to their offspring.
The case $p_0+p_2=1$ corresponds to a {\it binary tree}, the  constant average 
progeny case $p_0=p_2=1/2$ is called {\it critical}.

Shreve's calculations imply $T_{ij} = 2^{j-i-1}$, indicating, in particular,
that the side branching only depends on the difference between 
the branch orders, not on their absolute values.
Eiji Tokunaga \cite{Tok78} generalized this idea by introducing 
conditions  
\begin{subequations}
\label{eq:tok}
\begin{align}
\label{eq:tok1}
T_{ij}&=T_{j-i}\quad{\rm for~some~sequence~}\{T_k\ge 0\}_{k\ge 1},\\
\label{eq:tok2}
T_k&= a\,c^{k-1},\quad k\ge 1,~ a\ge0,~c>0.
\end{align}
\end{subequations}
The Tokunaga model \eqref{eq:tok} completely specifies a combinatorial tree 
shape (up to a permutation of side branch attachment within a given branch) with only two parameters $(a,c)$, hence suggesting a conventional modeling paradigm.
The empirical validity of the Tokunaga constraints \eqref{eq:tok} 
has been strongly confirmed for a variety of river networks at
different geographic locations \cite{Pec95,Tarb96,DR00,MTG10,ZZF13},
as well as in other types of data represented by trees, including 
botanical trees \cite{NTG97}, 
the veins of botanical leaves \cite{TPN98,PT00}, 
clusters of dynamically limited aggregation \cite{O92,NTG97}, 
percolation and forest-fire model clusters \cite{ZWG05,YNTG05}, 
earthquake aftershock sequences \cite{THR07,HTR08,Y13}, 
tree representation of symmetric random walks \cite{ZK12},
and hierarchical clustering \cite{GNT99}.
The conditions \eqref{eq:tok}, however, lack a theoretical justification.

This work suggests such a justification, by exploring a space $\cT$ of finite binary 
rooted trees with an arbitrary probability measure.
We introduce the {\it Horton pruning} -- cutting the tree leaves (Fig.~\ref{fig:pruning}), which is closely connected to the 
concept of Horton-Strahler orders. 
We showed in \cite{KZ16a} that prune invariant measures are abundant on $\cT$, and 
every prune invariant measure necessarily satisfies the constraint \eqref{eq:tok1}.
This motivates a {\it geometric branching process} (Def.~\ref{def:GBM}) that 
induces a prune invariant measure for an arbitrary sequence $\{T_k\}$.
A trajectory of this process is a random {\it geometric tree} (Def.~\ref{def:geom}).
The geometric tree that satisfies \eqref{eq:tok2} with 
$(a,c)=(1,2)$ is equivalent to the critical 
binary Galton-Watson tree.
A {\it time invariant} distribution preserves the average 
counts of branches of each order under 
the operation of a unit time shift in a geometric branching process (Def.~\ref{def:TI}, Fig.~\ref{fig:Depth}).
Our main result (Thm.~\ref{thm:main}) is that for the  
geometric trees, time invariance is equivalent to the Tokunaga constraint 
\eqref{eq:tok2} with $a=c-1$ and $c\ge1$.
Proofs are given in a separate section at the end of the paper.

\section{\label{sec_ss} Horton Prune Invariance}
Denote by $\cT$ the space of finite unlabeled rooted reduced binary trees
with no planar embedding, 
including the {\it empty tree} $\phi$ comprised 
of a single root vertex and no edges.
The existence of the root induces a parent-offspring relation for
each pair of adjacent vertices: the one closer to the root is called parent, 
the other -- offspring. 
The absence of labeling means that the tree vertices are not
distinguishable (other than by their position within a tree).
The absence of planar embedding refers to the lack of ordering (left/right)
among the offspring of the same parent.
The term {\it reduced} refers to the absence of vertices of degree 2.

\begin{figure}
\includegraphics[width=0.4\textwidth]{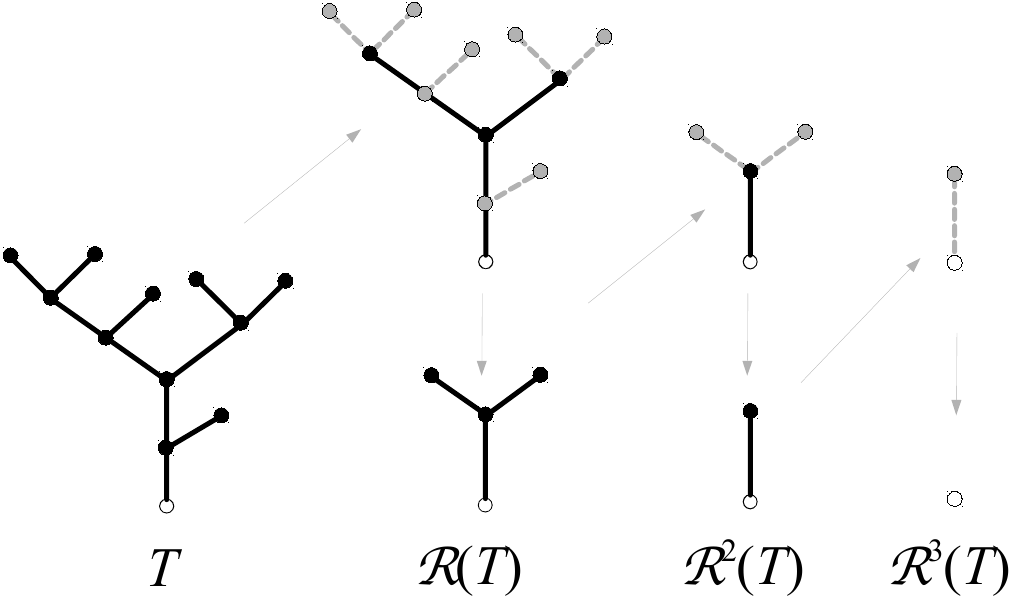}
\caption{\label{fig:pruning} Pruning: an illustration. 
The bottom row shows the trajectory $T\to\cR(T)\to\cR^2(T)\to\cR^3(T)=\phi$
of an initial tree $T$.
The top row highlights the edges and vertices (dashed gray) to be pruned 
at each step.
Notice two series reduction (elimination of a degree-2 vertex) that occur
after the first pruning.}
\end{figure}

The concept of tree self-similarity is related to the pruning operation
\cite{Pec95,BWW00,KZ16a}.
{\it Pruning} of a tree is a map $\cR:\cT\to\cT$,
whose value $\cR(T)$ for $T\ne\phi$ is obtained by removing
the leaves and their parental edges from $T$, followed by series reduction
(removing vertices of degree 2).
We also set $\cR(\phi)=\phi$. 
Pruning is an onto but not invertible operator,
with an infinite preimage $\cR^{-1}(T)$ for any $T\in\cT$.
Pruning induces a contracting flow on $\cT$. 
The trajectory of each tree $T$ under $\cR(\cdot)$ is uniquely
determined and finite (Fig.~\ref{fig:pruning}):
\[T\equiv\cR^0(T)\to \cR^1(T) \to\dots\to\cR^k(T)=\phi.\]

\begin{Def}[{{\bf Horton-Strahler orders}}]
\label{def:HS}
The Horton-Strahler order ${\sf ord}(T)\in\mathbb{Z}_+$ of a tree $T\in\cT$ is 
the minimal number of prunings necessary to erase the tree:
\be\label{eq:order}
{\sf ord}(T)=\min\left\{k\ge 0:\cR^k(T)=\phi\right\}.
\ee
\end{Def}

\begin{figure}
\includegraphics[width=0.35\textwidth]{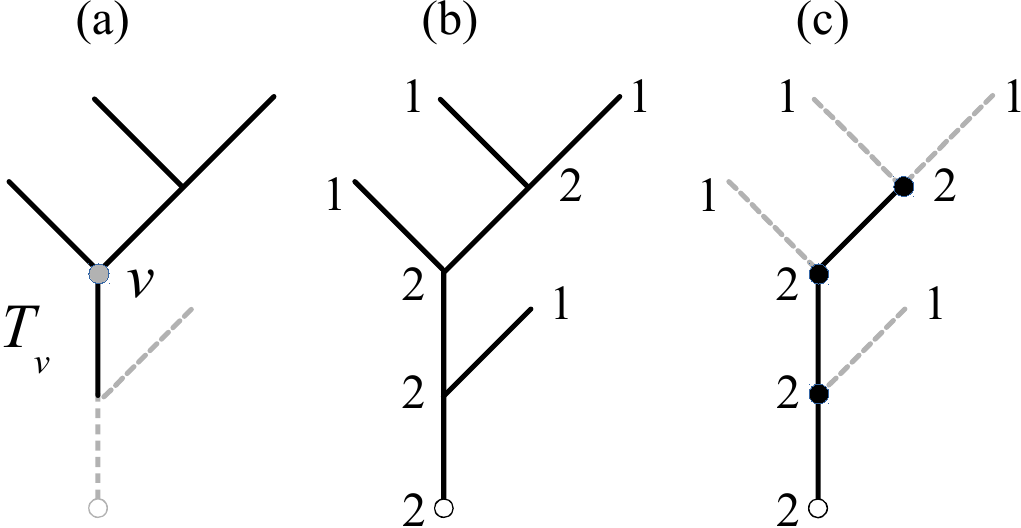}
\caption{\label{fig:Horton} Horton-Strahler ordering: an illustration. 
(a) Subtree $T_v$ (solid) descendant to vertex $v$ in $T$ (dashed gray).
(b) Horton-Strahler orders of tree vertices. 
The tree order ${\sf ord}(T)=2$ is assigned to the tree root.
(c) A branch of order ${\sf ord}=2$ consists of three vertices and their parental
edges (solid) in $T$ (dashed gray).}
\end{figure}

The order ${\sf ord}(v)$ of a vertex $v\in T$ is the order of the {\it descendant subtree} 
$T_v\subset T$ that includes all vertices descendant to $v$ (including $v$) in $T$, together with
their parental edges (Fig.~\ref{fig:Horton}a,b).
The order of the parental edge of $v$ is ${\sf ord}(v)$.
A {\it branch} is a collection of adjacent vertices/edges of the same order 
(Fig.~\ref{fig:Horton}c).
A {\it side branch} of order $\{ij\}$, $i<j$, is a branch of order $i$ that merges
a branch of a higher order $j$.
The tree in Fig.~\ref{fig:Horton}c has two side branches of order $\{12\}$; a situation
when two branches of the same order merge (and form a branch of a higher order) 
is {\it not} a side branch.

\begin{Def}[{{\bf Prune invariance}}]\label{def:prune}
Consider a probability measure $\mu$ on $\cT$ such that $\mu(\phi) = 0$.
Let $\nu(T)=\mu \circ \cR^{-1}(T) = \mu \big(\cR^{-1}(T)\big)$. (Note that $\nu(\phi)>0$.)
Measure $\mu$ is called {\it invariant} with
respect to the pruning operation if for any tree $T\in\cT$ we have 
\be
\label{def:pi}
\nu\left(T~|T\ne\phi\right)=\mu(T).
\ee
\end{Def}

We restrict our analysis to {\it coordinated} measures \cite{KZ16a}.
Informally, the coordination property asserts that all descendant 
subtrees of a given order $\kappa$ have the same probability structure, independent
of how deep or shallow they occur in a given tree. Formally,
a tree measure $\mu(\tau)$ is coordinated if, conditioned on ${\sf ord}(\tau)=K\ge 1$, 
for any $\kappa<K$, a uniformly selected descendant subtree of $\tau$ of order $\kappa$ is 
distributed as  $\mu(T|{\sf ord}(T)=\kappa)$. 
A tree measure $\mu(T)$, and a random tree drawn from $\mu(T)$,
is called {\it self-similar} if it is coordinated and prune invariant.

The authors have shown \cite{KZ16a} that for the coordinated measures, 
the prune invariance implies the order invariance constraint \eqref{eq:tok1}.
It was also shown that prune invariance implies a geometric
distribution of tree orders.
Accordingly, if one is interested in prune invariant trees, 
it is natural to consider measures specified by sequences $\{T_k\}$,
with a geometric order distribution.
The next section introduces a natural class of such measures.

\section{Geometric Branching Process}

We define here a multi-type discrete time branching process
that generates a prune invariant measure on $\cT$ that satisfy \eqref{eq:tok1}
for an arbitrary sequence $\{T_k\}$. 
In general, this can be done in numerous ways \cite{KZ16a}. 
We restrict ourselves to Markov processes that facilitate subsequent analysis.

We say that a random variable $X$ has a geometric distribution with parameter $r$, 
$X \stackrel{d}{=}{\sf Geom}(r)$,
if its probability mass function is given by
\[{\sf Prob}(X=k) = r(1-r)^k,\quad k=0,1,2,\dots\]

\begin{Def}[{{\bf Geometric branching process}}]
\label{def:GBM}
Consider a sequence of {\it Tokunaga coefficients} $\{T_k\ge0\}_{k\ge 1}$ and $0<p<1$.
Define $S_{K}:=1+T_1+\dots+T_{K}$ for $K\ge 0$ by assuming $T_0=0$.
The {\it geometric branching process} ${\cG}(s;T_k,p)\equiv\cG(s)$ describes a discrete time 
population growth:
\begin{itemize}
\item[(i)] The process starts at $s=0$ with a progenitor of order ${\sf ord}(\cG)$
such that ${\sf ord}(\cG)-1\stackrel{d}{=}{\sf Geom}(p)$.
\item[(ii)] At every time instant $s>0$, each population member of order
$K\in\{1,\dots,{\sf ord}(\cG)\}$ terminates 
with probability $q_K=S^{-1}_{K-1}$, independently of other members. 
At termination, a member of order $K>1$ produces two offspring of
order $(K-1)$; and a member of order $K=1$ terminates with leaving no offspring.
\item[(iii)] At every time instant $s>0$, each population member of order
$K\in\{1,\dots,{\sf ord}(\cG)\}$ survives (does not terminate) with probability $1-q_K=1-S^{-1}_{K-1}$,
independently of other members.
In this case, it produces a single offspring (side branch). 
The offspring order $i$, $1\le i<K$, is drawn from the distribution
\be
\label{eq:pj}
p_{K,i}=\frac{T_{K-i}}{T_1+\dots+T_{K-1}}.
\ee
\end{itemize}
\end{Def}

The geometric branching process is Markov
in the space of ornamented trees -- trees with vertex orders. 
The numbers of side branches
are independent for distinct branches.
Denote by $m_K$ the total number of side branches within a 
randomly selected branch of order $K$, and by $m_{K,i}$ the
number of side branches of order $i$, $i=1,\dots,K-1$, within
that branch.
Branches of order 1 have unit length and no side branches: $m_1=0$.
The property (ii) of the definition implies for any $K>1$  
\[ m_K\stackrel{d}{=}{\sf Geom}\left(S_{K-1}^{-1}\right),\]
with ${\sf E}[m_K]=S_{K-1}-1=T_1+\dots+T_{K-1}.$
Combining this with \eqref{eq:pj} we find 
\[ m_{K,i}\stackrel{d}{=}{\sf Geom}\left(\left[1+T_{K-i}\right]^{-1}\right),\]
with ${\sf E}\left[m_{K,i}\right] = T_{K-i}$.
The independence of branches implies
${\sf E}\left[N_{ij}\right] = {\sf E}\left[N_j\right]{\sf E}\left[m_{j,i}\right]$
and hence 
\[T_{ij} = \frac{{\sf E}[N_{ij}]}{{\sf E}[N_j]} 
= \frac{{\sf E}\left[N_j\right]{\sf E}\left[m_{j,i}\right]}{{\sf E}[N_j]}=T_{j-i}.\]
This means that a trajectory of the geometric branching process 
is a (random) tree ${\cG}(T_k,p)\in\cT$ that satisfies the condition \eqref{eq:tok1}
with the sequence $\{T_k\}$.
We notice that  essential elements of the resulting trees (tree order, branch lengths, 
numbers of side branches of a given order within a selected branch) are described by 
geometric laws -- hence the model name.

A tree ${\cG}(T_k,p)\in\cT$ generated by the geometric branching process can 
be equivalently defined via the following construction.

\begin{Def}[{{\bf Geometric tree}}]
\label{def:geom}
Consider a sequence $\{T_k\ge0\}_{k\ge 1}$ and $0<p<1$.
We say that a (random) tree ${\cG}(T_k,p)\in\cT$ is a {\it geometric} tree if and only if:
\begin{itemize}
\item[(i)] The tree order satisfies 
${\sf ord}(\cG)-1\stackrel{d}{=}{\sf Geom}(p)$.
\item[(ii)] The total number $m_K$ of side branches within a branch
of order $K\ge 2$ is ${\sf Geom}\left(S^{-1}_{K-1}\right)$.
\item[(iii)] For a branch of order $K\ge 2$, conditioned on $m_K=m$, 
the assignment of orders for $m$ side branches is done according to multinomial distribution
with $m$ trials and success probabilities $p_{K,j}$ of \eqref{eq:pj}
for any $j=1,\dots,K-1$.  
\item[(iv)] The number of side branches and their orders are independent 
in distinct branches.
\end{itemize} 
\end{Def}

The Def.~\eqref{def:geom} is convenient for recursive tree generation.
A tree of order $K=1$ consists of two vertices (root and leaf) connected by an edge.
To generate a random tree of general order, one first generates a random geometric
order $K$ and starts with a complete binary tree of depth $K$, which we call {\it skeleton}.
All leaves in the skeleton have the same depth $K$, and all vertices at depth $\kappa\le K$
have the same Horton-Strahler order $\kappa$.
Accordingly, the skeleton has order $K$ equal to the tree depth.
The final tree is obtained by adding side branches of lower orders to
every branch of the skeleton, according to the definition rules.
Every side branch is generated according to the same procedure.

It follows from the properties of a critical binary Galton-Watson
tree \cite{BWW00} that the geometric tree with $p=1/2$ and $T_k = 2^{k-1}$
is the critical binary Galton-Watson tree.
Notably, the next statement shows that the measure induced by a geometric tree
(and hence by a geometric branching process)
is prune invariant (Def.~\ref{def:prune}).
\begin{thm}[{{\bf Prune invariance}}]
\label{thm:PI}
Given an arbitrary sequence $\{T_k\ge0\}_{k\ge 1}$ 
and $0<p<1$,
the probability measure for the geometric random tree $\cG(T_k,p)$  
is invariant with respect to Horton pruning.
\end{thm}
\begin{proof}
See the Proofs section.
\end{proof}

A geometric tree measure is coordinated by construction.
The coordination together with prune invariance implies that the geometric 
measure is self-similar.

\begin{figure}
\includegraphics[width=0.4\textwidth]{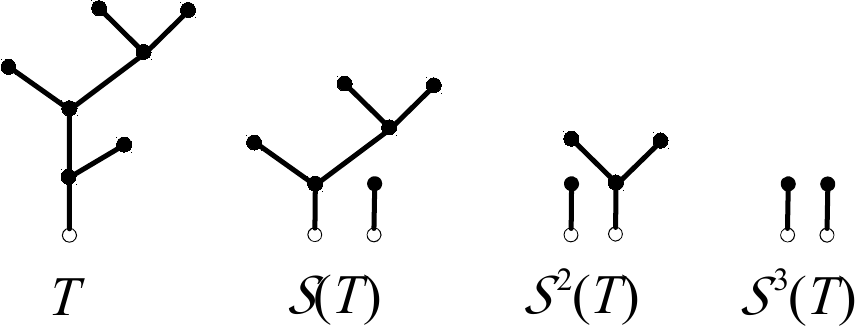}
\caption{\label{fig:Depth} Time shift: an illustration. 
The figure shows forests obtained by consecutive
application of the time shift operator $\cS$ to a tree $T$ shown on the left.
At every step, we remove the root edge from each existing tree. 
This terminates the trees of order ${\sf ord}=1$, and splits any other
tree in two new trees.}
\end{figure}

\section{Time Invariance}
Let $x_i(s)$, $i\ge 1$, denote the average number of vertices of order $i$
at time $s$ in the process $\cG(s)$, and
${\bf x}(s)=(x_1(s),x_2(s),\dots)^T$ be the state vector.
By definition we have
\[{\bf x}(0)=\pi:=\sum\limits_{K=1}^\infty p(1-p)^{K-1} {\bf e}_K,\]
where ${\bf e}_K$ are standard basis vectors. 
Furthermore, if $q_{a,b}$, $a\ge b$, denotes the probability that a vertex
of order
${\sf ord}=a+\mathbbm{1}_{\{a=b\}}$
that exists at time $s$ splits into a pair of vertices of 
orders $(a,b)$ at time $(s+1)$, then
\begin{eqnarray}
\label{eq:dyn1}
\lefteqn{x_K(s+1)= 2\,x_{K+1}(s)q_{K,K}}\nonumber\\
&+&x_K(s)(1-q_{K-1,K-1})
+\sum_{i=K+1}^{\infty}x_i(s)\,q_{i,K}.
\end{eqnarray}
The first term in the right-hand side of \eqref{eq:dyn1} corresponds to a split of an order-$(K+1)$
vertex into two vertices of order $K$, 
the second -- to a split of an order-$K$ vertex into a vertex
of order $K$ and a vertex of a smaller order, and
the third -- to a split of a vertex of order $i>K$ into 
a vertex of order $K$ and a vertex of order $i$. 
The geometric branching implies
\begin{eqnarray}
\label{joint}
q_{a,b}
=\left\{
\begin{array}{rl}
S_{a}^{-1} & \text{ if }a=b,\\
T_{a-b}S_{a-1}^{-1}& \text{ if }b<a.
\end{array}\right.
\end{eqnarray}

Summing up, the system \eqref{eq:dyn1} can be written as
\begin{eqnarray}
\label{eq:dyn2}
{\bf x}(s+1)-{\bf x}(s)&=&\mathbb{G}\mathbb{S}^{-1}{\bf x}(s),
\end{eqnarray}
where
\begin{equation} 
\label{eq:G}
\mathbb{G}:=\left[\begin{array}{ccccc}-1 & T_1+2 & T_2 & T_3 & \hdots  \\0 & -1 & T_1+2 & T_2 & \hdots  \\0 & 0 & -1 & T_1+2 & \ddots  \\0 & 0 & 0 & -1  & \ddots \\\vdots & \vdots & \ddots & \ddots  & \ddots \end{array}\right]
\end{equation}
and 
\[\mathbb{S} = {\sf diag}\{S_0,S_1,\dots\}.\]

In this setup, the {\it unit time shift} operator $\cS$, which advances the 
process time by unity, can be applied to 
individual trees and forests (collection of trees) -- it removes the root edge 
from a tree $T$.
A consecutive applications of $d$ time shifts to a tree $T$ is equivalent to removing
the vertices/edges at depth less than $d$ from the root (Fig.~\ref{fig:Depth}).
Next we define an invariance with respect to the time shift $\cS$.

\begin{Def}[{\bf{Time invariance}}]
\label{def:TI}
Process $\cG(s)$ is called {\it time invariant} if and only if
the state vector ${\bf x}(s)$ is invariant
with respect to a unit time shift $\cS$:
\be
\label{eq:TI}
{\bf x}(s) = {\bf x}(0)\equiv \pi ~~\forall s~\Longleftrightarrow~
\mathbb{G}\mathbb{S}^{-1}\pi = {\bf 0}.
\ee
\end{Def}

Now we formulate the main result of this paper.
\begin{thm}[{{\bf Critical Tokunaga process}}]
\label{thm:main}
A geometric branching process $\cG(s;T_k,p)$ is time invariant if and only if
\be
\label{eq:cTok}
p=1/2\text{  and } T_k=(c-1)c^{k-1} \text{ for any } c\ge 1.
\ee  
We call this family a critical Tokunaga process, and the respective trees -- critical
Tokunaga trees.
\end{thm}
\noindent Theorem~\ref{thm:main} is proven in the Proofs section
via solving a nonlinear system of equations that writes \eqref{eq:TI} in terms of ratios $S_k/S_{k+1}$.

Let $K$ be the order of a random tree $\cG(T_k,p)$, and, conditioned on  $K>1$, let $K_a,K_b$ be the orders 
of its two subtrees, $T_a$ and $T_b$, rooted at the internal vertex closest to the root, randomly and uniformly permuted.
We call $T_a$ and $T_b$ the {\it principal} subtrees of $T$.

\begin{cor}
\label{cor1}
Let $\cG$ be a critical Tokunaga tree. 
Then the distribution of the principal subtree $T_a$ 
(and hence $T_b$) matches that of the initial tree $\cG$.
The distributions of $T_a$ and $T_b$ are independent if and only if 
$c=2$.
\end{cor}
\begin{proof}
See the Proofs section.
\end{proof}

\section{Discussion}
We suggest a measure-theoretic justification for the algebraic
Tokunaga constraint \eqref{eq:tok} with $a=c-1$ and $c\ge 1$, 
in terms of invariance of geometric branching process $\cG(s;T_k,p)$ 
with respect to a time shift. 
As a result, there appears a family of critical Tokunaga branching processes
\eqref{eq:cTok} characterized by invariance to both the Horton pruning (erasing
a tree from leaves to the root) 
and time shift (erasing a tree from the root to leaves).

\paragraph{Extending properties of critical Galton-Watson tree.} 
Our results extend some remarkable properties of the critical binary 
Galton-Watson process ($T_k=2^{k-1}$) to a one-parametric family 
of critical Tokunaga processes ($T_k=(c-1)c^{k-1}$):

(i) Burd et al.~\cite{BWW00} have shown that condition \eqref{eq:tok1} characterizes 
the invariance with respect to the Horton pruning in the family of binary 
Galton-Watson trees. 
Our work \cite{KZ16a} established that
\eqref{eq:tok1} follows
from the prune invariance for any coordinated measure
(recall that any Galton-Watson measure is coordinated, but
not vice-versa).
Moreover, coordination is a necessary property -- one can 
construct a multitude of uncoordinated prune invariant measures 
that do not satisfy the constraint \eqref{eq:tok1}.

(ii) Among binary Galton-Watson trees, the  
condition \eqref{eq:tok1}, and hence prune invariance,
only holds in the {\it critical} case ($p_0=p_2=1/2$).
In this case \eqref{eq:tok2} is also satisfied with $a=c-1=1$ \cite{BWW00}.
We show here that the condition \eqref{eq:tok2} with $a=c-1$ is equivalent to time 
invariance (which implies criticality) in a broader family of geometric trees.

(iii) The forest obtained from descendant subtrees rooted at every vertex of 
a critical binary Galton-Watson tree $T$ approximates the frequency 
structure of the entire space of trees, as the order of $T$ increases.
The same property holds for the trees generated by a critical Tokunaga 
branching process with any $c\ge 1$.

\paragraph{Independence of Horton-Strahler ordering.}
Our results relate the conditions \eqref{eq:tok1} and
\eqref{eq:tok2}, expressed in terms of the Horton-Strahler orders, 
to a discrete erasure of a tree from the leaves (Horton pruning)
or from the root (time invariance), respectively. 
Such erasures, that do not use the notion of order,
provide more natural process-generating constraints.

\paragraph{Ease of simulation.}
Generation of geometric trees for arbitrary parameters $(p,\{T_k\})$ is
easily implemented on a computer (see discussion after Def.~\ref{def:geom}).
This model may hence facilitate analysis in a range of simulation-heavy problems,
from structure and transport on river networks to phylogenetic tree analysis.

\paragraph{Horton law.} 
The Horton law -- a geometric decay of the branch counts $N_i$ 
-- is among the most fundamental regularities found in observed trees,
both static and dynamic
\cite{Pec95,Tarb96,NTG97,RIR01,DR00,KZ17,ZFG10}.
It was shown in \cite{KZ16a} that a strong version of the Horton law
is satisfied in a self-similar tree, provided 
$\limsup_{k\to\infty} T_k^{1/k}<\infty$.
In particular, for a critical Tokunaga tree with 
$T_k=(c-1)c^{k-1}$ and $c\ge 1$ the Horton ratios $N_i/N_{i+1}$
asymptotically converge to $R_b=2c$, for any $i\ge 1$, as the 
tree order increases.

\paragraph{Fractal dimension of trees.}
Recall that a fractal dimension of a tree is defined
as $d = \ln{R_b}/\ln{R_r}$ \cite{NTG97}, with the asymptotic Horton ratios
$R_b=\lim N_i/N_{i+1}$ and $R_r=\lim r_{i+1}/r_i$,
where $r_i$ is the average length of a branch of order $i$.
If we assume that the lengths of edges in a geometric tree are 
independent equally distributed random variables with a unit mean, then ${\sf E}(r_i)=S_{i-1}$. 
For the critical Tokunaga model this gives $R_r = c$ and 
the respective fractal dimension 
\be
\label{eq:fd}
d_c=\frac{\ln{2c}}{\ln{c}}=1+\ln(2)/\ln(c).
\ee
Hence, $1\le d_c< \infty$ depending
on the value of $c$. 
In particular, $d_{2^{1/k}}=1+k$, which corresponds to a volume-filling
tree in $(1+k)$ dimensions for any $k\ge 1$.

\paragraph{Horton ratio.} 
The Horton ratio reported in the observed large trees is approximately within $2.8<R_b<6$.
The expression \eqref{eq:fd} for the fractal dimension of a critical Tokunaga tree,
together with $R_b=2c$, gives
$\log_2(R_b) = {\frac{d_c}{d_c-1}}.$
Accordingly, the observed empirical range of the Horton ratios is reproduced
by the critical Tokunaga trees with $1.4< c < 3$ and dimensions $1.6 < d_c < 3$.
This covers all the tree dimensions that may exist in a 3-dimensional world, excluding
the range $d_c<1.6$ that corresponds to almost ``linear'', and 
probably less studied, trees. 

\paragraph{On the condition $a=c-1$.}
The constraint \eqref{eq:tok2} with $a=c-1$ that defines the critical Tokunaga 
process appears in the Random Self-similar Network model \cite{VG00,MTG10}, 
introduced in hydrological studies.
At the same time, the multiple documented instances when the Tokunaga constraint
\eqref{eq:tok2} is tested in observations and models
typically refer to $a\ne c-1$ \cite{Pec95,Tarb96,DR00,MTG10,ZZF13,NTG97,O92,GNT99,ZWG05}. 
Although the observed deviations are not large in absolute value
(e.g., $c\approx 2.6$, $a\approx 1.1$ for river basins \cite{Pec95,ZZF13}),
they seem to be statistically significant \cite{ZZF13}.
It remains an open problem to either find a theoretical justification for
the condition \eqref{eq:tok2} with $a\ne c-1$, or otherwise explain the apparent 
deviation from the time invariance that leads to $a\ne c-1$ in observed
trees.

\section{Proofs}
\subsection{Proof of Theorem~\ref{thm:PI}}
\begin{proof}
Consider a random geometric tree $\cG\equiv\cG(T_k,p)$.
Let $T:=\cR(\cG)$ be the respective pruned tree.
We show below that tree $T$ is geometric,
i.e., it satisfies properties (i)-(iv) of Def.~\eqref{def:geom}.
Observe that the Horton pruning $\cR$ eliminates branches of order 1 in the initial tree $\cG$ 
and decreases the order of every non-leaf branch by unity. 
Hence, the event $\{{\sf ord}(\cG)>1\}$ is the same as
$\{T\ne\phi\}$.

(i) We have
\[{\sf P}\left({\sf ord}(T)=k|T\ne\phi\right)= 
\frac{{\sf P}\left({\sf ord}(\cG)=k+1\right)}{{\sf P}\left({\sf ord}(\cG)>1\right)},\]
which implies $\{{\sf ord}(T)|T\ne\phi\}-1\stackrel{d}{=}{\sf Geom}(p)$.

(ii) Pruning $\cR$ acts as a {\it Bernoulli thinning} on the number $m_K$ of side branches within 
every branch of order $K\ge 2$, with thinning probability $p_{K,1}$ of \eqref{eq:pj}.
Recall that application of a Bernoulli thinning with removal probability $1-q$ to 
a geometric random variable with parameter $r$ results in a geometric random variable with parameter
\[r_{\rm thinned} = \frac{r}{q(1-r)+r}.\]
To obtain the distribution for the number $m_K$ of side branches in a branch of 
order $K$ of a pruned tree $T$, we apply a Bernoulli thinning to a branch of order $(K+1)$
in the original tree $\cG$.
Using $r=S_{K}^{-1}$ and $1-q=p_{K+1,1}$ shows that $m_K \stackrel{d}{=}{\sf Geom}(S_{K-1}^{-1})$,
which establishes (ii).

(iii) Follows from the properties of multinomial distribution.

(iv) Follows from independence of the branch structure in $\cG$. 
\end{proof}

\subsection{Proof of Theorem~\ref{thm:main}}
\noindent A proof of Theorem~\ref{thm:main} is based on the following 
two lemmas proven in separate sections below.

\begin{lem}
\label{lem:system}
A geometric branching process $\cG(s)$ is time invariant if and only if
$p=1/2$ and the sequence $\{T_k\}$ solves the following (nonlinear) system
of equations:
\be
\label{eq:system}
\frac{S_0}{S_k}=\sum_{i=1}^{\infty} 2^{-i}\frac{S_i}{S_{k+i}}\quad\text{for all }k\ge 1.
\ee
\end{lem}

\noindent Let $a_k=S_k/S_{k+1}\le 1$ for all $k\ge 0$.
Then, for any $i\ge 0$ and any $k>0$ we have $S_i/S_{k+i} = a_i\,a_{i+1}\dots a_{i+k-1}$.
The system \eqref{eq:system} rewrites in terms of $a_i$ as
\begin{align*}
{1 \over 2}a_1+{1 \over 4}a_2+{1 \over 8}a_3+\hdots &=a_0, \nonumber \\
{1 \over 2}a_1a_2+{1 \over 4}a_2a_3+{1 \over 8}a_3a_4+\hdots &=a_0a_1, \nonumber \\
{1 \over 2}a_1a_2a_3+{1 \over 4}a_2a_3a_4+{1 \over 8}a_3a_4a_5+\hdots &=a_0a_1a_2,\nonumber  
\end{align*}
and so on, which can be summarized  as
\begin{equation}\label{eq:a}
\sum\limits_{j=1}^\infty {1 \over 2^j}\prod\limits_{k=j}^{n+j-1}a_k =  \prod\limits_{k=0}^{n-1}a_k,
\text{ for all } n \in \mathbb{N}.
\end{equation}
Assume that $\{a_0,a_1,a_2,\dots\}$ is a solution to system \eqref{eq:a}.
Then $\{1,a_1/a_0,a_2/a_0,\dots\}$ is also a solution, since each equation only includes multinomial terms of the same degree.
\begin{lem}
\label{lem:NL}
The system \eqref{eq:a} with $a_0=1$ has a unique solution
$a_0=a_1=a_2=\hdots=1.$
\end{lem}
\noindent Lemma~\ref{lem:NL} implies
$a_k = S_k/S_{k+1} = 1/c$ for some $c\ge 1$.
Hence
$S_1 = 1+T_1 = c$ and $T_1 = c-1.$
Furthermore, 
\[S_{k+1} = c\,S_k =c^k\]
and, accordingly,
\[T_{k+1} = S_{k+1}-S_k = (c-1)c^{k-1},\]
which completes the proof.

\subsection{Proof of Lemma~\ref{lem:system}}
\begin{proof}
Assume that the process is time invariant. 
Then the process progeny
is constant in time and equals unity:
\[\norm{\pi}_1=\sum_{k=1}^{\infty}p(1-p)^{k-1} = 1.\]
Observe that in one time step, every vertex of order ${\sf ord}=1$
terminates, and any vertex of order ${\sf ord}>1$ splits in two.
Hence, the process progeny at $s=1$ is
\[2\sum_{k=2}^{\infty}p(1-p)^{k-1} = 2(1-p)=1,\]
which implies $p=1/2$.
Accordingly, $p(1-p)^{k-1}=2^{-k}$ and
the time invariance \eqref{eq:TI} takes the following coordinate form
\begin{eqnarray}
\label{eq:cdn}
-\frac{2^{-k}}{S_{k-1}}&+&2^{-(k+1)}\frac{T_1+2}{S_k}\nonumber\\
&+&\sum_{i=k+2}^{\infty}2^{-i}\frac{T_{i-k}}{S_{i-1}}=0,\text{ for all } k\ge 1.
\end{eqnarray}
Multiplying \eqref{eq:cdn} by $2^k$ and observing that $T_k=S_k-S_{k-1}$ 
we obtain:
\begin{eqnarray}
\label{eq:10}
&\displaystyle-\frac{1}{S_{k-1}}&+\frac{1}{2}\frac{T_1+2}{S_k}
+\sum_{i=2}^{\infty}2^{-i}\frac{T_{i}}{S_{k+i-1}}=0,\nonumber\\
&\displaystyle\frac{1}{S_{k-1}}&-\sum_{i=1}^{\infty}2^{-i}\frac{S_i}{S_{k+i-1}}\nonumber\\
&=&\frac{1}{S_k}-\frac{1}{2\,S_k} -\sum_{i=2}^{\infty}2^{-i}\frac{S_{i-1}}{S_{k+i-1}},\nonumber\\
&\displaystyle\frac{1}{S_{k-1}}&-\sum_{i=1}^{\infty}2^{-i}\frac{S_i}{S_{k+i-1}}\nonumber\\
&=&\frac{1}{2}\left(\frac{1}{S_k}-\sum_{i=1}^{\infty}2^{-i}\frac{S_i}{S_{k+i}}\right),
\end{eqnarray}
We prove \eqref{eq:system} by induction.
Base: For $k=1$ we have
\begin{eqnarray*}
\frac{1}{2}&=&\frac{1}{2\,S_1}+\sum_{i=1}^{\infty}2^{-(i+1)}\frac{S_i-S_{i-1}}{S_i},\\
1&=&\frac{1}{S_1}+\sum_{i=1}^{\infty}2^{-i}-\sum_{i=1}^{\infty}2^{-i}\frac{S_{i-1}}{S_i},
\end{eqnarray*}
which gives
\[\frac{1}{S_1}=\sum_{i=1}^{\infty}2^{-i}\frac{S_i}{S_{i+1}}.\]
Step: Assume that the statement is proven for $(k-1)$.
Then the left-hand side of \eqref{eq:10} vanishes, 
and the right-hand part rewrites as \eqref{eq:system}.
This establishes necessity.

Conversely, we showed that the system \eqref{eq:system} 
is equivalent to \eqref{eq:TI} in case $p=1/2$.
This establishes sufficiency.
\end{proof}

\subsection{Proof of Lemma~\ref{lem:NL}}

\begin{proof} We consider two cases. 
\medskip
\noindent

{\it Case I.} Suppose the sequence $\{a_j\}$ has a maximum:
there exists an index $i \in \mathbb{N}$ such that 
$a_i =\max\limits_{j \in \mathbb{N}}a_j.$
Define
\[w_{j,\ell}:= {1 \over 2^j}\prod\limits_{k=j}^{\ell+j-1}a_k \left[\prod\limits_{k=0}^{\ell-1}a_k\right]^{-1}.\]
Using $n=\ell$ in (\ref{eq:a}) we obtain that for any $\ell \in \mathbb{N}$,
\begin{equation}\label{eq:L}
\sum\limits_{j=1}^\infty w_{j,\ell}  =  1,
\end{equation}
and using $n=\ell+1$ we find that an arbitrary $a_\ell$ is the weighted average of $\{a_{\ell+j}\}_{j=1,2,\hdots}$:
\begin{equation}\label{eq:L+1}
\sum\limits_{j=1}^\infty w_{j,\ell}\,a_{\ell+j} =  a_\ell.
\end{equation}
Hence, since $a_i =\max\limits_{j \in \mathbb{N}}a_j$,
$$a_i=a_{i+1}=a_{i+2}=a_{i+3}=\hdots =a.$$
Similarly, letting $\ell=i-1$ in (\ref{eq:L}) and (\ref{eq:L+1}), we obtain $a_{i-1}=a$. Recursively, by plugging in $\ell=i-2, ~i-3, \hdots$, we show that
$$a_1=a_2=\hdots =a_{i-1}=a_i=a_{i+1}=\hdots =a.$$
Finally, ${1 \over 2}a_1+{1 \over 4}a_2+{1 \over 8}a_3+\hdots=1$ implies $a=1$.

\medskip
\noindent
{\it Case II.} Suppose there is no $\max\limits_{j \in \mathbb{N}}a_j$. Let 
$U:=\limsup\limits_{j \rightarrow \infty}a_j.$
From (\ref{eq:a}) we know via cancelation that
\begin{eqnarray}
\label{eq:U2}
{1 \over 2}a_n&+&{1 \over 4}{a_n a_{n+1} \over a_1}+{1 \over 8}{a_n a_{n+1}a_{n+2} \over a_1a_2}+\hdots\nonumber\\
&+&{1 \over 2^{n-1}}{ \prod\limits_{k=n}^{2n-2}a_k \over \prod\limits_{k=0}^{n-2}a_k}
+\sum\limits_{j=n}^\infty {1 \over 2^j}{ \prod\limits_{k=j}^{n+j-1}a_k \over \prod\limits_{k=0}^{n-1}a_k}=  1.
\end{eqnarray}
Thus, $2^{-1}\,a_n <1$ and $U \leq 2$.
The absence of maximum implies $a_j <U \leq 2$ for all $j \in \mathbb{N}$.

\medskip
\noindent
Plugging $n+1$ in (\ref{eq:a}), we obtain
\begin{eqnarray*}\label{eq:Un+1}
\lefteqn{\left({1 \over 2}a_n\right) a_{n+1}+\left({1 \over 4}{a_n a_{n+1} \over a_1}\right)a_{n+2}+\hdots}\\
&+&{1 \over 2^{n-1}}{ \prod\limits_{k=n}^{2n-2}a_k \over \prod\limits_{k=0}^{n-2}a_k}a_{n+j-1}+\sum\limits_{j=n}^\infty {1 \over 2^j}{ \prod\limits_{k=j}^{n+j-1}a_k \over \prod\limits_{k=0}^{n-1}a_k}a_{n+j}=  a_n.
\end{eqnarray*}

Thus, since $a_j <U$ for all $j \in \mathbb{N}$,
\begin{eqnarray*}
\lefteqn{\left({1 \over 2}a_n\right) a_{n+1}+\left({1 \over 4}{a_n a_{n+1} \over a_1}\right)U+\hdots}\\
&+&{1 \over 2^{n-1}}{ \prod\limits_{k=n}^{2n-2}a_k \over \prod\limits_{k=0}^{n-2}a_k}U+\sum\limits_{j=n}^\infty {1 \over 2^j}{ \prod\limits_{k=j}^{n+j-1}a_k \over \prod\limits_{k=0}^{n-1}a_k}U > a_n
\end{eqnarray*}
which simplifies via (\ref{eq:U2}) to
\begin{equation}\label{eq:anan+1}
\left({a_n \over 2}\right) a_{n+1}+\left(1-{a_n \over 2}\right)U  > a_n.
\end{equation}
For all $\varepsilon \in (0,1)$, there are infinitely many $n \in \mathbb{N}$ such that $a_n >(1-\varepsilon)U$.
Then, for any such $n$, the above inequality (\ref{eq:anan+1}) implies
$$a_{n+1}>2-{2 \over a_n}U+U >2-{2\varepsilon \over 1-\varepsilon}+U =\big(1-\varphi(\varepsilon)\big)U,$$
where 
$$\varphi(x):={2x \over (1-x)U}.$$
Let $\varphi^{(k)}=\varphi \circ \hdots \circ \varphi$. 
Repeating the argument for any given number of iterations $K \in \mathbb{N}$, we obtain
\begin{eqnarray*}
a_{n+2}&>&\big(1-\varphi^{(2)}(\varepsilon)\big)U, \quad a_{n+3}>\big(1-\varphi^{(3)}(\varepsilon)\big)U,\\ 
&&\hdots, a_{n+K}>\big(1-\varphi^{(K)}(\varepsilon)\big)U.
\end{eqnarray*}
Thus, given any $K \in \mathbb{N}$, fix $\varepsilon \in (0,1)$ small enough so that such that $\varphi^{(k)}(\varepsilon)\in (0,1)$ for all $k=1,2,\hdots,K$. Then, taking $n>K$ such that $a_n >(1-\varepsilon)U$, we obtain from (\ref{eq:U2}) that
\begin{eqnarray*}
1 &>& {1 \over 2}a_n+{1 \over 4}{a_n a_{n+1} \over a_1}+{1 \over 8}{a_n a_{n+1}a_{n+2} \over a_1a_2}+\hdots\\
&+&{1 \over 2^{K+1}}{ \prod\limits_{k=n}^{n+K}a_k \over \prod\limits_{k=0}^Ka_k}\\
&>& {1 \over 2}(1-\varepsilon)U+{1 \over 4}{(1-\varepsilon)\big(1-\varphi(\varepsilon)\big)U^2 \over U}+\hdots\\
&+&{1 \over 2^{K+1}}{ U^{K+1}\prod\limits_{k=0}^K\big(1-\varphi^{(k)}(\varepsilon)\big) \over U^K}.
\end{eqnarray*}
Now, since $\varepsilon$ can be chosen arbitrarily small,
$$1 \geq \left(1-{1 \over 2^{K+1}}\right)U.$$
Finally, since $K$ can be selected arbitrarily large, we have proven that $1\geq U$. However, this will contradict the assumption of {\it Case II}. Indeed, if
$a_j <U \leq 1$ for all $j \in \mathbb{N}$, then
$${1 \over 2}a_1+{1 \over 4}a_2+{1 \over 8}a_3+\hdots <1,$$
contradicting the first equation in the statement of the theorem. Thus, the assumptions  of {\it Case II} cannot be satisfied.  
We conclude that there exists a maximal element in the sequence $\{a_j \}_{j=1,2,\hdots}$ as assumed in {\it Case I}, implying the statement of the theorem.
\end{proof}

\subsection{Proof of Corollary~\ref{cor1}}
\begin{proof}
Let ${\sf ord}(\cG)$ denote the (random) order of a random geometric tree $\cG$. 
We apply a unit time shift $\cS$ to the tree $\cG$.
Conditioned on ${\sf ord}(\cG)>1$, at instant $s=1$ there exist exactly two 
vertices that are the roots of the principal subtrees $T_a$ and $T_b$. 
Since the trees $T_a$ and $T_b$ have the same distribution, their roots
have the same order distribution.
Denote by $y_k$ the probability that the tree $T_a$ has order $k\ge 1$
and let ${\bf y}=(y_1,y_2,\dots)^T$.
By Thm.~\ref{thm:main} criticality is equivalent to time invariance.
In a critical tree we have $p=1/2$, which, together with time invariance,
implies
\[{\bf x}(0)={\bf x}(1)=2{\bf y}\pi_1+{\bf 0}(1-\pi_1)
=2{\bf y}\frac{1}{2}+{\bf 0}\frac{1}{2}={\bf y}.\]
This establishes the first statement.

The second statement follows from examining the joint distribution $q_{a,b}$
of \eqref{joint}.
\end{proof}

\section*{Acknowledgement}
The authors are grateful to the Editor and two referees whose comments helped to improve
the initial version of this work.
The authors acknowledge financial support from the National Science Foundation, 
awards DMS 1412557 (Y.K.) and  EAR 1723033 (I.Z.)



\end{document}